\documentclass[reqno]{amsart}
 \usepackage{amsfonts,amssymb,eucal}
\usepackage{amsmath} 
\usepackage{amsfonts}
\usepackage{amssymb}
\usepackage{amscd}
\usepackage{amsthm}

\def\ee{\mathrm{e}}

\numberwithin{equation}{section}

\input{cyracc.def}

\numberwithin{equation}{section}
\newcommand{\Ad}{\ensuremath{{\mbox{\rm{Ad}}}}}% Ad- representation
\newcommand{\tr}{\ensuremath{{\mbox{\rm{Tr}}}}}% Ad- representation
\newcommand{\ad}{{\mbox{\rm ad}}}
 
\def\i{\mathrm{i}}

\newcommand{\mb}[1]{{\mbox{\boldmath{$#1$}}}}% mathematical bold
\newcommand{\mc}[1]{{\mathcal{#1}}}% mathematical caligraphic
\newcommand{\got}[1]{{\mathfrak{#1}}}% gothic with mbox for  mathematic
\newcommand{\db}[1]{{\mathbb{#1}}}% double

\newcommand{\pa}{\partial}

\newcommand{\R}{\ensuremath{\mathbb{R}}}
\newcommand{\C}{\ensuremath{\mathbb{C}}}

\newcommand{\T}{\ensuremath{\text{SU}(3)/S(\text{U}(1)\times
\text{U}(1)
 \times \text{U}(1))}}
\newcommand{\Hi}{\ensuremath{\mathcal{H}}}% Hilbert space
\newcommand{\Z}{\ensuremath{\mathbb{Z}}}
\newcommand{\fl}{\ensuremath{{\mathcal{F}}_{\Hi}}}%Bargmann H space
\newtheorem{lemma}{Lemma}[section]
\newtheorem{corollary}[lemma]{Corollary}

\newtheorem{Proposition}[lemma]{Proposition}

\newtheorem{Remark}[lemma]{Remark}
\newtheorem{Comment}[lemma]{Comment}
\newcommand{\K}{K\"ahler~}

\title[A holomorphic 
representation of  the Jacobi algebra]{A holomorphic 
representation of the multidimensional  Jacobi  algebra}
\author{Stefan  Berceanu}
\address[Stefan  Berceanu]{National
 Institute for Physics and Nuclear Engineering\\
         Department of Theoretical Physics\\
         PO BOX MG-6, Bucharest-Magurele, Romania}
\email{Berceanu@theor1.theory.nipne.ro}


\begin{document}
\begin{abstract}We present  a holomorphic  representation of the
Jacobi algebra
$\mathfrak{h}_n\rtimes \mathfrak{sp}(n,\R )$ by first order
differential 
operators with polynomial coefficients on the manifold
$\mathbb{C}^n\times  
\mathcal{D}_n$. We construct  the Hilbert space of holomorphic functions
on which these differential  operators act.
\end{abstract}

\maketitle
%\noindent
%\tableofcontents
{\small \tableofcontents}
%\newpage
%\noindent
%\tableofcontents
%\newpage
\section*{Introduction}

The {\it coherent states} (CS), invented by Schr\"odinger
in the early days of quantum
mechanics, offer a bridge between quantum and classical physics, see
 e.g.  \cite{perG,zh1,ali} and  references therein.
On the other hand, the group-theoretic generalization \cite{per}   from 
the standard Glauber's field CSs  (attached to the
 Heisenberg-Weyl (HW) group) to  coherent states associated  to
any Lie group, attracted mathematicians, as a
convenient
way to  representation theory and  quantization  \cite{berezin2}. Today,
one speaks about CS-type groups \cite{lis,neeb}, a
 class of groups which
contains all compact groups, all simple hermitian groups, certain
solvable groups and also mixed groups such as the semidirect product of
the HW group  and the symplectic group. The latter,
 called after Eichler and Zagier {\it Jacobi group}  \cite{ez},
 captured 
the interest of mathematicians  working in   representation theory
(e.g. Satake
\cite{sat}, Berndt and Schmidt \cite{bs}
and Neeb \cite{neeb}),
 and  automorphic forms (e.g.  K\"ahler in his last three papers
\cite{cal1,cal2,cal3}).
Independently,  the   Jacobi group has already  been studied by
physicists    under the name of {\it Schr\"odinger} or {\it Hagen group}
 in the seventies \cite{nied,hagen}. On the other hand,
 the Jacobi group 
  governs the  so called {\it squeezed states} in Quantum Optics. 
 I shall point out in several places below
    these parallel contexts
in which the Jacobi group appears. 

The aim of this paper  is to study holomorphic representations of the
Jacobi group on  a certain homogeneous  K\"ahler manifold
$M$ attached
to the  group. Propositions \ref{final} and \ref{finalf}  express
our main results.
 Strictly speaking, mathematicians have produced
general schemes \cite{sat,neeb}
 from which, some of  the formulas presented here should  be
obtained as particular cases. 
From this point of view, the paper addresses mainly to a  reader  more
familiar with the methods of the Theoretical or Mathematical
Physics, who does not want to read hundreds of pages to know what
 the CS-groups really are, but who needs explicit formulas, like those
presented  by our Remark \ref{difa}. Such  formulas are
 useful  when  studying
the equations of motion on $M$ generated by  Hamiltonians which are
algebraic functions  in
the generators of the  group
\cite{berezin1,sbcag,sbl}, here the Jacobi group.  However, while
the scalar product or the reproducing kernel (\ref{KHK}), related to 
the K\"ahler
potential
on $M$, is known \cite{neeb}, as far as we know, the resolution
 of unity 
(\ref{ofi}) has not
been written down explicitly for this case.
 The  formulas furnished by the
Propositions \ref{final}, \ref{finalf} can be used for Berezin's 
quantization
\cite{berezin2}. Remark \ref{re9}  shows the equivalence of the 
present approach
with the similar theorems obtained in reference \cite{ez}.   

The starting point of the story is generally   considered to be
  the Segal-Bargmann-Fock realization \cite{bar}
$a\mapsto \frac{\partial}{\partial z};~~ a^+\mapsto z$ of the
canonical commutation relations (CCR) $[a,a^+]=1$  on the
symmetric Fock space $\fl :=\Gamma^{\mbox{\rm{hol}}}(\C,\frac{\i}{2
\pi}\exp(-|z|^2)dz\wedge d\bar{z})$  attached to the Hilbert space
$\Hi :=L^2(\R, dx)$. In fact, much earlier Sophus Lie
 had done the
differential  realization of the generators
$K^{0,-,+}$ of the group  $\text{SU}(1,1)$,  
$K^{-}\mapsto \frac{\pa }{\pa w}$; $K^{0}\mapsto k+ w\frac{\pa }{\pa w}$;
$K^{+}\mapsto 2kw + w^2 \frac{\pa }{\pa w}$ 
on the unit disk  $w\in \mc{D}_1=\text{SU}(1,1)/\text{U}(1)$. 
 In \cite{holl} we have presented
an explicit   realization
of a holomorphic representation  for the  the CS-group  $G^J_1$,
which is the semidirect product of the real three dimensional HW group
 $H_1$ with the group $\text{SU}(1,1)$, called the Jacobi group
\cite{ez}.  
  These formulas contain both  the standard   
  Segal-Bargmann-Fock representation \cite{bar}
 and the  realization
of the generators of the group $\text{SU}(1,1)$ already mentioned.  
 
 At this
point, let us mention that the present investigation is part of a
larger program, we started earlier.  
  In reference \cite{last} we 
advanced the conjecture that  the generators
of all CS-groups \cite{lis,neeb} admit representations by 
first order differential operators with holomorphic
polynomials coefficients on CS-manifolds. Let us recall some of
our own earlier
progress in this field.
 Our method \cite{morse}
 permits to get the holomorphic differential action of the generators
of a continuous unitary representation $\pi$
   of a Lie  group $G$ with the Lie algebra $\got{g}$  on a homogeneous
space $M=G/H$. Following Perelomov \cite{per},
we consider 
homogeneous manifolds $M$ realized as K\"ahler 
CS-orbits \cite{sbctim,last,sinaia}.   Previously,  we
 produced explicit representations 
for hermitian groups, using CS
 based on compact \cite{sbcag}
 and noncompact \cite{sbl}
hermitian symmetric spaces, and also on \K CS-orbits of semisimple Lie
groups \cite{sbctim,sinaia}. In all such situations the  differential
action of the generators of  the group $G$ on holomorphic
functions defined on   a \K homogeneous orbit 
  $G/H$ can be written down as a sum of two terms, the first one
as a polynomial $P$, and, the second one, as sum of partial derivatives times
some polynomials $Q$-s. On hermitian symmetric spaces the degrees of the
polynomials $P$ and $Q$ are less than  or equal to two \cite{sbcag,sbl}.
 We have analyzed the simplest example of a  nonsymmetric
homogeneous manifold, \T, where the maximum degree of the polynomials
 $P, Q$   already equals three
\cite{sbctim,sinaia}.

More precisely, the 
present paper is devoted to a concrete realization of a differential
holomorphic representation of the Lie algebra
 $\got{g}^J_n:= \got{h}_n\rtimes
\got{sp}(n,\R )$
on the  homogeneous space $M :=  H_n/\R\times \text{Sp}(n,\R
)/\text{U}(n)$, where
$H_n$ denotes the $(2n+1)$-dimensional HW group. The complex 
$n(n+3)/2$-dimensional  manifold
 $M $  is realized as  $\mc{D}^J_n:=\C^n\times\mc{D}_n$, where
$\mc{D}_n$
is the $n(n+1)/2$-dimensional Siegel ball, endowed with a \K structure
deduced  from the scalar product of two CSs based on $M$.

Now we digress a little about the connection between the
 present paper and other fields. As it had   already been mentioned, 
the physical object associated with the Jacobi group
 are the  squeezed states of the Quantum Optics
 \cite{zh1,siv,dod,dr}, discovered  already in the
early  days of Quantum Mechanics \cite{ken}.  
   Instead of starting from  a
matrix representation    of the Jacobi group (see e.g. p. 182 in
\cite{kir1}), in our
presentation we use methods similar to those of the squeezed states
\cite{stol}.
It is well known that for the
harmonic oscillator CSs the uncertainties in momentum and position
are equal to $1/\sqrt{2}$ (in units of $\hbar$).
 ``The squeezed
states'' \cite{ken,stol,lu,bi,yu,ho,wa}
are the states for which the uncertainty in  position is less than
$1/\sqrt{2}$. The squeezed states are a particular class of ``minimum
uncertainty states'' (MUS) \cite{mo}, 
i.e. states which saturate the Heisenberg
uncertainty relation. In the present paper we do not
insist on  possible ``physical'' 
applications \cite{dr} of our paper to  the squeezed states.
 Let us just mention that ``Gaussian pure states''
(``Gaussons'') \cite{si} are more general MUSs. In fact, as it
 was shown in
\cite{ali}, these states are CSs based on the manifold
$\mc{X}^J_n:= \mc{H}_n\times\R^{2n}$, where $\mc{H}_n$ is the Siegel
upper  half
plane  $\mc{H}_n:=\{Z\in M(n,\C)| Z=U+\i V, U,V\in M(n,\R), \Im (V)>0, U^t=U;
V^t=V\}$.
$M(n,R)$ denotes the $n\times n$ matrices with entries in $R$,
$R=\R~\text{or}~\C$, and $X^t$ denotes the transpose of the matrix
$X$. In  \cite{sbj}
 we have started the generalization of CSs attached to the
 Jacobi group $G^J_1$
 to the (``multidimensional'') Jacobi group
$G^J_n$. The connection of  our
construction of coherent states based on $\mc{D}^J_n$
\cite{sbj}  with
the Gaussons of \cite{si} is a subtle one and it should be investigated
separately. $\mc{D}_n$ denotes the Siegel ball  $\mc{D}_n:= \{Z\in
M(n,\C)|Z=Z^t, ~ 1-Z\bar{Z}> 0 \}$. In \S \ref{unu1} we  
recall  the clue  of this connection
  in the  case
$n=1$ \cite{holl}, 
 which is offered by the K\"ahler-Berndt's construction
\cite{cal1,cal2,cal3,cal,bern,bs}.

Furthermore we try to make a technical presentation of 
the content of the paper.
 Firstly, the  notation referring to
the HW group is fixed in \S \ref{HWF}. In \S \ref{simpl}
some known facts are recalled: the definition of the symplectic
algebra and the symplectic
group, the Gauss and Cartan decomposition, the differential action of
the generators, the scalar product. \S \ref{JTR} is devoted to the
Jacobi group. We introduce coherent states associated to the Jacobi group 
$G^J_n$
 based on the homogeneous space $\mc{D}^J_n$. The explicit
formulas giving the differential action of the generators of the
Jacobi group are given in Remark \ref{difa}. Remark \ref{hoho},
which can be considered as generalizing the
Holstein-Primakoff-Bogoliubov type equations, implies in Remark
\ref{remark6} the action of the group $\text{Sp}(n,\R )$ on $H_n/\R
\equiv\C^n$. 
 Lemma \ref{lema6} is very important for our construction: it connects
the normalized ``squeezed'' vector \cite{stol}
$\Psi_{\alpha,W}=D(\alpha)
S(W)e_0$, $\alpha\in\C^n$, $W\in\mc{D}_n$ with Perelomov's 
un-normalized vector $e_{z,W}$.
As a consequence, we get in Remark \ref{comm1} the expression of the
scalar product of two CS vectors associated to the Jacobi group, based
on the manifold $\mc{D}^J_n$. As we had already  emphasized, this
expression is already known (see p. 532 in
\cite{neeb}, or the article \cite{hno}
and (5.28) in \cite{sat}), but here we  present a simple proof. 
Proposition \ref{mm1} allows to find the action of the Jacobi group
$G^J_n$ on the 
base manifold $\mc{D}^J_n$ and the composition law in
the Jacobi group. Remark  \ref{re9} identifies our expression with the
one obtained in 
context of Jacobi forms \cite{ez}. In
\S \ref{two} we find out  the K\"ahler two-form
$\omega$ on the manifold $\mc{D}^J_n$.
Applying  the technique
of Chapter IV from Hua's book \cite{hua} and a lemma from Berezin's
paper \cite{berezin2},  we determine the Liouville measure on 
 $\mc{D}^J_n$. This permits the explicit
 construction of  the symmetric Fock space
$\mc{F}_K$ 
attached to the reproducing kernel $K$, summarized in Proposition
\ref{final}. Proposition \ref{finalf} gives
the continuous unitary holomorphic representation $\pi_K$ of the group
$G^J_n$ on $\mc{F}_K$.  The \S \ref{unu1} recalls some of
the further results established for the Jacobi group $G^J_1$. Finally, 
we discuss
 the connection between
 different contexts in which the same group appears under
the names of Jacobi,
Schr\"odinger or Hagen.  

\section{The Heisenberg-Weyl group}\label{HWF}

The Heisenberg-Weyl   group  $H_n$ is the nilpotent  group with the 
$2n+1$-dimensio\-nal real  Lie algebra isomorphic to the 
algebra 
 \begin{equation}\label{nr0}\got{h}_n =
<\i s 1+\sum_{i=1}^n (x_i a_i^+-\bar{x}_ia_i)>_{s\in\R ,x_i\in\C}
,\end{equation} 
 where ${ a}_i^+$ (${ a}_i$) are  the boson creation
(respectively, annihilation)
operators which verify the CCR 
\begin{equation}\label{baza1}
[a_i,a^+_j]=\delta_{ij}; ~ [a_i,a_j] = [a_i^+,a_j^+]= 0 . 
\end{equation}
The vacuum verifies the relations:
\begin{equation}\label{vacuma}
a_ie_o= 0, i=1, \cdots, n.
\end{equation}

The displacement operator
\begin{equation}\label{deplasare}
D(\alpha ):=\exp (\alpha a^+-\bar{\alpha}a)=\exp(-\frac{1}{2}|\alpha
|^2)  \exp (\alpha a^+)\exp(-\bar{\alpha}a),
\end{equation}
verifies the composition rule:
\begin{equation}\label{thetah}
D(\alpha_2)D(\alpha_1)=\ee^{\i\theta_h(\alpha_2,\alpha_1)}
D(\alpha_2+\alpha_1) , 
~\theta_h(\alpha_2,\alpha_1):=\Im (\alpha_2\bar{\alpha_1}) .
\end{equation}
  Here we have used the notation
 $\alpha \beta= \sum_i \alpha_i\beta_i$, where $\alpha = (\alpha_i )$.
The composition law of the HW group $H_n$ is:
\begin{equation}\label{clh}
(\alpha_2,t_2)\circ (\alpha_1,t_1)=(\alpha_2+\alpha_1, t_2+t_1+
\Im (\alpha_2\bar{\alpha_1})). 
\end{equation}
If  we identify $\R^{2n}$ with $\C^n$, $(p,q)\mapsto \alpha$:
\begin{equation}
\label{change}
\alpha = p +\i q, ~p,q\in\R^n,
\end{equation} then $$ \Im (\alpha_2\bar{\alpha_1})=
(p^t_1, q^t_1)J\left(\begin{array}{c}p_2 \\ q_2\end{array}\right),
\text{where} ~
J=\left(\begin{array}{cc} 0 & 1\\ -1 & 0\end{array}\right) .$$

\section{The symplectic group}\label{simpl}
\subsection{The symplectic algebra}
The real symplectic Lie algebra $\got{sp}(n, \R )$ is a real form of
the simple Lie algebra $ \got{sp}(n, \C )$ of type $\got{c}_n$
and $ X\in \got{sp}(n, \R )$ if
 $$  X^tJ+JX=0 ~\text{or}~
X=\left( \begin{array}{cc} a & b \\c & -a^t\end{array}\right), ~b=b^t,
c=c^t,$$
where $a,b,c\in M(n,\R )$, and similarly for $ \got{sp}(n, \C )$.

In the complex realization   (\ref{change}), to $X\in
\got{sp}(n,\R )$,  $X\in M(2n,\R)$ corresponds $X_{\C}\in \got{sp}(n,\C
)\cap \got{u} (n,n)$
\begin{equation}
X_{\C}=\left(\begin{array}{cc} a & b \\ \bar{b} &
\bar{a}\end{array}\right),
\end{equation}
where
\begin{equation}\label{cond1} a^*=-a, ~ b^t= b
\end{equation} (cf. theorems in
\cite{itzik}, \cite{bar70}, \cite{fol}). So, we consider 
the  realization of the Lie algebra of the group $\text{Sp}(n,\R )$:
\begin{equation}\label{nr1}
\got{sp}(n,\R)=
<\sum_{i,j=1}^n (2a_{ij}K^0_{ij} +
b_{ij}K^+_{ij}-\bar{b}_{ij}K^-_{ij})>, \end{equation}
where the matrices $a=(a_{ij})$, $b=(b_{ij})$ verify the conditions
(\ref{cond1}).

The  generators $K^{0,+,-}$ verify the  commutation
relations 
\begin{subequations}\label{baza2}
\begin{eqnarray}
 [K_{ij}^-,K_{kl}^-] & = & [K_{ij}^+,K_{kl}^+]=0 \label{baza21}, \\
 2[K_{ij}^-,K_{kl}^+] & = & K^0_{kj}\delta_{li}+
K^0_{lj}\delta_{ki}+K^0_{ki}\delta_{lj}+K^0_{li}\delta_{kj}
\label{baza22}, \\
2[K^-_{ij},K^0_{kl}] & = & K_{il}^-\delta_{kj}+K^-_{jl}\delta_{ki}
 \label{baza23}, \\
 2[K^+_{ij},K^0_{kl}] & = & -K^+_{ik}\delta_{jl}-K^+_{jk}\delta_{li}
 \label{baza24}, \\
 2[K^0_{ji},K^0_{kl}] & = & K^0_{jl}\delta_{ki}-K^0_{ki}\delta_{lj}
\label{baza25}.  
\end{eqnarray}
\end{subequations}
%We  consider the matrix realization
%\begin{equation}\label{nr2}
%K_{ij}^0=\frac{1}{2}\left(\begin{array}{cc} e _{ij}&0 \\0 &-e_{ji}\end{array}\right),~
%K_{ij}^+=\frac{i}{2}\left(\begin{array}{cc} 0 & e _{ij}+e _{ji} \\ 0 &0 
%\end{array}\right)~,
%K _{ij}^-=\frac{i}{2}\left(\begin{array}{cc}0 &0 \\  e _{ij}+e _{ji}
% & 0 \end{array}\right) .
%\end{equation}

With the notation: $\mb{X}:= d\pi (X)$, we have the correspondence:
$\got{sp}(n,\R )\ni X \mapsto \mb{X}$, where 
the real symplectic  Lie algebra $\got{sp}(n, \R )$ is
 realized as  $\got{sp}(n,\C)\cap \got{u}(n,n)$% as  in 
%eqs. (\ref{nr1}), (\ref{cond1}), (\ref{nr2}):

\begin{equation}\label{algebra}
X= \left(\begin{array}{cc} a & b \\ \bar{b} & \bar{a} \end{array} \right)
\leftrightarrow  
  ~~~\mb{X} =  \sum_{i,j=1}^n (2a_{ij}\mb{K}^0_{ij} +
z_{ij}\mb{K}^+_{ij}-\bar{z}_{ij}\mb{K}^-_{ij}) ,~ b=\i z,
\end{equation}
where (\ref{cond1}) is verified.
%\newpage
In Table 1 we give the realization of the generators of the real
symplectic group in matrices, as operators obtained via the derived
representation, and a bi-fermion realization. 
\begin{table}
\caption[Table I]{{\it The generators of the symplectic group: operators, matrices,
 and bifermion operators}}
\begin{center}
\begin{tabular}{||c|c|c||}\hline
~$\mb{K}^+_{ij}$ ~&
~ $ K_{ij}^+=\frac{\i}{2}\left(\begin{array}{cc} 0 & e _{ij}+e _{ji} \\ 0 &
0 
\end{array}\right)~$~ &
~ $\frac{1}{2}a^+_ia^+_j$~\\ \hline
~$\mb{K}^-_{ij}$ ~ & ~
$K _{ij}^-=\frac{\i}{2}\left(\begin{array}{cc}0 &0 \\  e _{ij}+e _{ji}
 & 0 \end{array}\right) $ & ~$ \frac{1}{2}a_ia_j $~\\ \hline
~$\mb{K}^0_{ij}$~ &~ $K_{ij}^0=\frac{1}{2}\left(\begin{array}{cc} e _{ij}&0 \\0 &
-e_{ji}\end{array}\right) $~ & ~$\frac{1}{4}(a^+_ia_j+a_ja^+_i)$~\\ \hline
\end{tabular}
\end{center}
\end{table}
\subsection{The symplectic group}
 For $g\in \text{GL}(2n,\R )$, we have
\begin{equation}\label{spr}
g\in \text{Sp}(n,\R )~\leftrightarrow ~ g^tJg=J. 
\end{equation}
 If in  (\ref{spr})
 $g\in \text{GL}(2n,\C )$, then $g\in \text{Sp}(n,\C )$. We
remind also that $g\in \text{U}(n,n)$ iff $gKg^*=K$, where
$K=\left(\begin{array}{cc}
1 & 0 \\ 0 & -1\end{array}\right)$.

Under the identification (\ref{change}) of $\R^{2n}$ with $\C^n$, we
have the correspondence  
$$A\in M(2n,\R )\rightarrow A_{\C}\in M(2n, \R )_{\C},
~A_{\C}=WAW^{-1},~ W=\frac{1}{\sqrt{2}}\left(\begin{array}{cc} 1 & \i 1
\\
1 & -\i 1\end{array}\right),$$ where
$$M(2n, \R )_{\C}= \left\{\left(\begin{array}{cc} P & Q \\ \bar{Q}
&\bar{P} \end{array}\right), P, Q\in M(n,\C )\right\}.$$
We extract  from \cite{bar70}, \cite{fol}
\begin{Remark} To every $g\in {\rm{Sp}}
(n,\R )$ as in {\rm{(\ref{spr})}},
 $g\mapsto g_c\in \rm{Sp}(n, \C)\cap
\rm{U}(n,n)$, or denoted just $g$ 
 \begin{equation}\label{dg}
g= \left( \begin{array}{cc}a & b\\ \bar{b} &
\bar{a}\end{array}\right),
\end{equation}
where
\begin{subequations}\label{simplectic}
\begin{eqnarray}
aa^*- bb^* & = & 1;~ ab^t=ba^t, \\
a^*a-b^t\bar{b} & = & 1;~ a^t\bar{b}=b^*a . 
\end{eqnarray}
\end{subequations}
\end{Remark}
\begin{Remark}
The linear canonical transformations, i.e. the transformations
 which leaves invariant
 {\rm{(\ref{baza1})}}, are given by elements of the group 
$\rm{Sp}(n,\R )$
under the realization $\rm{Sp}(n, \C )\cap \rm{U}(n,n)$.
\end{Remark}
See also Remark \ref{hoho} below.

If  $g\in \text{Sp}(n,\R)$ is given by  (\ref{dg}),
then
\begin{equation}\label{dg1}
g^{-1}= \left( \begin{array}{cc}a^* & -b^t\\ -b^* & a^t
\end{array}\right) .
\end{equation}

{\bf Gauss decomposition.}
Let us consider an element $g\in \text{Sp}(n,\R)$. The following relations are true: 
\begin{equation}\label{dg11}
 g = \left(\begin{array}{cc} a & b\\ \bar{b} & \bar{a}
\end{array}\right) = 
\left(\begin{array}{cc} 1 & Y\\ 0  & 1 \end{array}\right) 
\left(\begin{array}{cc} \gamma &  0\\ 0  & \delta \end{array}\right) 
\left(\begin{array}{cc} 1 & 0 \\ Y' & 1 \end{array}\right) ,
\end{equation}
or
\begin{equation}\label{dg12}
 g = \left(\begin{array}{cc} a & b\\ \bar{b} & \bar{a}
\end{array}\right) = 
\left(\begin{array}{cc} 1 & 0\\ U  & 1 \end{array}\right) 
\left(\begin{array}{cc} \gamma' &  0\\ 0  & \delta' \end{array}\right) 
\left(\begin{array}{cc} 1 & U' \\ 0 & 1 \end{array}\right) ,
\end{equation}
where
\begin{equation}
Y= b\bar{a}^{-1}; ~ Y'= \bar{a}^{-1}\bar{b};\delta = \bar{a};
~ \gamma =
a-b\bar{a}^{-1}\bar{b}=(a^*)^{-1}=(\delta ^t)^{-1};  
\end{equation}
\begin{equation}
1-YY^*=(aa^*)^{-1}>0; Y=Y^t; 1-Y'Y'^*=((a^*a)^t))^{-1}>0;  Y'=Y'^t; 
\end{equation}
\begin{equation}
U=\bar{Y}= \bar{b}a^{-1}; U'=\bar{Y}'=a^{-1}b;
\gamma'=\bar{\delta}=a;\delta'= \bar{\gamma} = (a^t)^{-1}.
\end{equation}

{\bf Cartan decomposition}.
Let us consider an element $g\in \text{Sp}(n,\R)$. The following relations are true: 
\begin{equation}\label{notatie1}
 g =  \left(\begin{array}{cc} a & b\\ \bar{b} & \bar{a}
\end{array}\right) = \left(\begin{array}{cc} m & n\\ p & q\end{array}\right)
 \left(\begin{array}{cc} v & 0\\ 0 & \bar{v}\end{array}\right)=
\exp  \left(\begin{array}{cc} 0 & Z\\ \bar{Z} & 0\end{array}\right) 
\left(\begin{array}{cc} v & 0\\ 0 & \bar{v}\end{array}\right),
\end{equation}
where
\begin{equation}\label{notatien}
m=\cosh \sqrt{Z\bar{Z}}; ~
 n= \frac{\sinh \sqrt{Z\bar{Z}}}{\sqrt{Z\bar{Z}}}Z = 
Z\frac{\sinh \sqrt{\bar{Z}Z}}{\sqrt{\bar{Z}Z}};~ p=\bar{n};~ q=\bar{m},
\end{equation}
or
\begin{equation}
m = (1-YY^+)^{-1/2};  ~ n= (1-YY^+)^{-1/2}Y; ~
 v= (1-YY^+)^{1/2}.
\end{equation}
$Z$ and $Y$ above are related by formulas (\ref{u3}) and 
(\ref{u4}) below, with the correspondence $Z\leftrightarrow Z$,
$Y\leftrightarrow W$, where $Y=b\bar{a}^{-1}$, i.e. 
\begin{equation}\label{u6}
Z  =   \frac{\mbox{\rm{arctanh}}\sqrt{YY^+}}{\sqrt{ YY^+}}Y = 
\frac{1}{2\sqrt{YY^+}}\log\frac{1+ \sqrt{ YY^+}}{1- \sqrt{ YY^+}};
 ~Y=b\bar{a}^{-1}.
\end{equation}

{\bf Hermitian symmetric spaces}. 
We briefly recall some well known facts about hermitian symmetric
spaces \cite{helg,wolf}. We use the notation:\\
\noindent $X_n$- hermitian symmetric space of noncompact type, 
$X_n= G_0/K$;\\
$X_c$- compact dual of $X_n$, $X_c= G_c/K$; \\
$G_0$- real hermitian group;\\ $G=G_0^{\C}$ - 
the   complexification of $G_0$;\\  $P$ -
 a parabolic subgroup of $G$;\\ $K$ - maximal compact subgroup of $G_0$;\\
 $G_c$ - compact real form of $G$;\\  
  The compact manifold $X_c$ of
$\frac{n(n+1)}{2}$-complex dimension has a 
complex structure inherited
from the identification $X_c= G_c/K = G/P $.\\ The group $G_c$ acts
transitively on $X_c$  with isotropy group $
K= G_0\cap P= G_c\cap P$.\\
  $X_n
=G_0/K= G_n(x_0)$ is open in $X_c$, where $x_0$ is a base point of $G$
corresponding to $K$.\\  $X_c$ includes $X_n$ under Borel
embedding $X_n\subset X_c$: $gK\rightarrow gP, g\in G_0$. 

In our case: $G_0= \text{Sp}(n,\R )$, $G=\text{Sp} (n,\C )$,
 $G_c= \text{Sp} (n)= \text{Sp} (n, \C )\cap
\text{U}(2n)\subset \text{SU}(2n)$,
 $K= \text{U}(n)$, and $$P=\left\{\left(\begin{array}{cc}
a & 0\\ c & d
\end{array}\right);~ a^tc=c^ta, ~ a^t d = 1 \right\}.$$
Let
\begin{equation}
\got{m}^+= \left\{\left(\begin{array}{cc}0 & b \\ 0 & 0
\end{array}\right), ~ 
b^t=b 
\right\} .
\end{equation}
Then
\begin{equation}
Z\rightarrow \hat{Z}=\left(\begin{array}{cc}0 & Z\\ 0 & 0\end{array}
\right), 
~\xi (Z)=( \exp\hat{Z})x_0, 
\end{equation}
and 
$\xi$ maps the symmetric $n\times n$ matrices $Z$, $1-Z\bar{Z}>0$
 of $\got{m}^+$ onto a dense open subset of $X_c$ that contains 
$X_n$. This gives the Harish-Chandra embedding: $X_n\subset \xi
(\got{m}^+)\subset
X_c$. 
The non-compact  hermitian  symmetric 
space $X_n= \text{Sp}(n, \R )/\text{U}(n)$ admits a
realization  as a bounded homogeneous domain,  the Siegel ball  $\mc{D}_n$
\begin{equation}\label{dn}
\mc{D}_n:=\left\{W\in  M (n, \C )| W=W^t, 1-W\bar{W} > 0\right\}.
\end{equation}
$X_n$ is a hermitian symmetric space of type CI (cf. Table V. p. 518 in
\cite{helg}), identified with the  symmetric bounded domain of type
II, 
$\got{R}_{II}$
in Hua's notation \cite{hua}.

\subsection{Coherent states for the symplectic group}

Coherent states associated to the real symplectic group were
considered in several  references as 
\cite{mon,berezin1}, see also \S 8 in \cite{perG}. 
We consider a particular case of the positive discrete series
representation \cite{knapp} of $\text{Sp}(n,\R )$. The vacuum is chosen such that
\begin{subequations}\label{vac}
\begin{eqnarray}
\mb{K}^+_{ij} e_0 & \not= & 0 ,\label{vac1}\\
\mb{K}^-_{ij} e_0 & = & 0 ,\label{vac2}\\
\mb{K}^0_{ij} e_0 & = & \frac{k}{4}\delta_{ij} e_0. \label{vac3} 
\end{eqnarray}
\end{subequations}
We have the relations:
\begin{equation}
\pi \left(\begin{array}{cc} v & 0\\ 0 & \bar{v}\end{array}\right) 
e_0= (\det v)^{k/2}e _0, ~v\in \text{U}(n);
\end{equation}

\begin{equation}
\pi \left(\begin{array}{cc} a & 0\\ 0 & d \end{array}\right) e_0 =
 (\det a)^{k/2} e_0, ~ da^t =1, 
~ v\in \text{U}(n);
\end{equation}

\begin{equation}
\pi \left(\begin{array}{cc} a & 0\\ 0 & d \end{array}\right) = 
\exp (\sum 2A_{ij}\mb{K}^0_{ij}),~ a = \exp A .
\end{equation}

We introduce some notation, and we find out:
\begin{subequations}
\begin{eqnarray}
\underline{S}(Z) &  = &\exp (\sum z_{ij}
\mb{K}^+_{ij}-\bar{z}_{ij}\mb{K}^-_{ij}),~ Z = (z_{ij}); \label{u1}\\
\label{u2}
 S(W) & = &\exp (W\mb{K}^+)\exp (\eta \mb{K}^0) \exp (-\bar{W}K^-); \\
& = &
\exp (-\bar{W}\mb{K}^-)\exp (-\eta \mb{K}^0)\exp (W\mb{K}^+);\\
 W & = & Z\tanh \frac{\sqrt{Z^*Z}}{\sqrt{Z^*Z}}\label{u3};\\
Z  & = & \frac{\mbox{\rm{arctanh}}\sqrt{WW^*}}{\sqrt{WW^*}}W = 
\frac{1}{2}\frac{1}{\sqrt{WW^*}}
\log \frac{1 + \sqrt{WW^*}}{1-\sqrt{WW^*}}\label{u4};\\
\eta & = & \log (1-WW^*) = -2 \log \cosh \sqrt {ZZ^*}\label{u5}. 
\end{eqnarray}
\end{subequations}
We have $\underline{S}(Z)= S(W)$.  
In (\ref{u6}) $Y$ is that from the Gauss decomposition (\ref{dg11}).

Perelomov's un-normalized CS-vectors are:  
\begin{equation}\label{ccs}
e_Z:= \exp ({\sum z_{ij}\mb{K}^+_{ij}}) e_0 = \pi 
 \left(\begin{array}{cc} 1 & \i Z \\ 0 &  1 \end{array}\right) e_0 ,~ 
Z=(z_{ij}), ~ Z=Z^t.
\end{equation}

Let us consider an element $g\in \text{Sp}(n,\R )$. 
\begin{Remark}The following relations between the normalized and 
un-normalized Perelomov's CS-vectors hold:
\begin{equation}\label{r1}
\underline{S}(Z) e_0 = \det (1-WW^*)^{k/4}
e_{W},  % ~ W=Z\frac{\tanh\sqrt {Z^*Z} }{\sqrt {Z^*Z}},
\end{equation}
\begin{equation}\label{r2}
e_g:\!= \!\pi (g) e_0\!=\!
\pi \!\left(\!\begin{array}{cc} \!a \!& \!b\!\\ \!c\! &
\!d\!\end{array}\!\right)\! e_0 \! = \!  
(\det \!\bar{a})^{-k/2}\! e_{Z}\!= \!
 \left( \frac{\det \!a}{\det
{\!\bar{a}}}\right)^{\frac{k}{4}}\!\underline{S}(Z)e_0,
 Z\!=\!-\i b\bar{a}^{-\!1}, 
\end{equation}
\begin{equation}\label{r3}
S(g)e_{W/\i}= \det (Wb^*+a^*)^{-k/2}e_{Y/\i},
\end{equation}
where $W\in\mc{D}_n$, and  $Z\in\C^n$
 in {\em  (\ref{r2})} are related by equations  {\em   (\ref{u3})},
{\em (\ref{u4})},  and
the linear-fractional action of the group $\rm{Sp}(n,\R )$ on
the unit ball $\mc{D}_n$ in {\em   (\ref{r3})} is
\begin{equation}\label{r4}
Y:= g\cdot W 
=(a \, W+ b)(\bar{b}\,W  +\bar{a})^{-1} = (Wb^*+a^*)^{-1}(b^t+Wa^t).
\end{equation}
\end{Remark}
Using the results of \cite{bgras,cluj}, it can be proved
\begin{Remark}If $\underline{S}(Z)$ is defined by {\rm{(\ref{u1})}},
 then: 
\begin{subequations}
\begin{eqnarray}
\underline{S}(Z_2)\underline{S}(Z_1) & = & \underline{S}(Z_3)e^{2A\mb{K}^0}
\label{unuu1};\\
\underline{S}(Z_2)\underline{S}(Z_1)e_0 & = & 
(\det v)^{k/2}\underline{S}(Z_3) e_0
\label{unuu3};
\end{eqnarray}
\end{subequations}
\begin{subequations}
\begin{eqnarray}
W_3 & = & \!(1\!-\!W_1W_1^*\!)^{-1/2}\!(W_1\!+\!W_2\!)\!(1\!+\!W^*_1W_2\!)^{-1}
\!(1\!-\!W_1^*W_1\!)^{1/2};
\label{unuu2}\\
\ee^{A} & = & ~v ~= (MM^+)^{-1/2}M; \label{unuu5} \\
M & = & (1-W_1W_1^*)^{-1/2}(1+W_1W_1^*)(1-W_2W_2^*)^{-1/2};\label{thetas}\\
\det (v) & = &
\left[ \frac{\det (1+W_1W^*_2)}{\overline{\det} (1+W_1W^*_2)}\right]^{1/2},
\label{unuu4}
\end{eqnarray}
\end{subequations}
where in  {\em  (\ref{unuu2})} $W_i$ and $Z_i$, $i=1,2,3$,
 are related by the relations {\em (\ref{u3}), (\ref{u4})}.
\end{Remark}

\subsection{The differential action for the group $\text{Sp} (n,\R)$}
 
We consider CS-vectors given by (\ref{ccs})
$$e_W=\ee^{\mb{X}}e_0, 
~\mb{X}=\sum w_{ij}\mb{K}^+_{ij}; ~W=(w_{ij}), ~ W=W^t .$$
It is easy to see that:
\begin{subequations}
\begin{eqnarray}
\mb{K}^+_{kl} e_W & = & \frac{\pa}{w_{kl}}e _W, \\
\mb{K}^0_{kl} e_W & = & \left(\frac{k}{4}\delta_{kl}+
w_{il}\frac{\pa}{\pa w_{ik}}\right)e_W, \\
\mb{K}^-_{kl} e_W & = &
\left(\frac{k}{2}w_{kl}+
w_{al}w_{ik}\frac{\pa}{\pa w_{ai}}\right) e_W.
\end{eqnarray}
\end{subequations}
The proof is based on the general formula
\begin{equation}\label{for}
\Ad (\exp X) =\exp (\ad_X),
\end{equation}
valid  for Lie algebras $\got{g}$, which here  we write down explicitly as
\begin{equation}\label{bazap}
A\ee^X=\ee^X(A-[X,A]+\frac{1}{2}[X,[X,A]]+\cdots ) .
\end{equation}
The differential action for the group  $\text{Sp} (n,\R)$ is obtained:

\begin{subequations}\label{astae}
\begin{eqnarray}
 \db{K}^- & = & \frac{\pa}{\pa W}, \\
 \db{K}^+ & = & \frac{k}{2}W + W \frac{\pa}{\pa W} W, \\
 \db{K}^0 & = & \frac{k}{4}1 + \frac{\pa}{\pa W} W.
\end{eqnarray}
\end{subequations}
 We have used the convention:
$$\left[\left(\frac{\pa}{\pa W}W\right) f(W)\right]_{kl}:=
\frac{\pa f(W)}{\pa w_{ki}}w_{il}, ~ W= (w_{ij}).$$
Formulas of the type (\ref{astae}) have been obtained by
 several authors \cite{krpa,rowe}.   
  
\subsection{The scalar product}
The scalar product of two coherent state vectors indexed by the points
of the Siegel ball $\mc{D}_n$ is:
\begin{equation}\label{kzz}
K(Z,Z')=(e_Z,e_{Z'})_{\Hi}= \det (1- Z'Z^+)^{-k/2}.
\end{equation}
The K\"ahler two-form on $\mc{D}_n$ is:
\begin{equation}\label{kel} 
\omega   =  \i \frac{k}{2}\tr [(1-W\bar{W})^{-1}d W
\wedge (1-\bar{W}W)^{-1} d\bar{W}].
\end{equation}
The Jacobian $J=\det dW'/dW $
of the transformation $W'=g\cdot W = f(g,W)$ such that
$W'(W_1)=0$, modulo an irrelevant unitary transformation, is obtained
(see \cite{bgras})
from 
 $$dW'=(1 -W_1\bar{W}_1)^{-\frac{1}{2}}dW(1-\bar{W}_1W_1)^{-\frac{1}{2}}. $$
The density $Q=|J|^2$ of the $\text{Sp}(n,\R )$-invariant volume  is $$ Q = \det
(1-W\bar{W})^{-(n+1)}, $$
and the Bergman
 kernel, modulo the volume of the domain
$\mc{D}_n$ (cf. \cite{hua}) is
 $$K_{\mc{D}_n}(X,Y)=\det(1-X\bar{Y})^{-(n+1)}.$$
The scalar product in the space of functions   
$f_{\psi}(z)= (e_{\bar{z}},\psi )\in\mc{F}_{\Hi}$ is
\begin{equation}\label{ofiW}
(\phi ,\psi )_{\fl}= \Lambda_1 \! \int_{1-W\bar{W}>0}\!
\bar{f}_{\phi}(W)f_{\psi}(W) \det (1-W\bar{W})^q dW, ~q=\frac{k}{2}- n-1, 
\end{equation}
where $\Lambda_1 =J_n^{-1}(q)$, with $J_n(p)$ ($p>-1$) defined by:
\begin{equation}\label{JJJ}
J_n(p) = \frac{\pi^{\frac{n(n+1)}{2}} }{(p+1)\cdots (p+n)}\cdot
\frac{\Gamma (2p+3)
\Gamma (2p+5)\cdots \Gamma (2p+ 2n -1)}{\Gamma (2p+ n+ 2)
\Gamma (2p + n+3)\cdots \Gamma (2p+ 2n)}.
\end{equation}
We can write  (\ref{JJJ}) in the  form
\begin{equation}\label{JJJ1}
J_n(p)=2^n\pi^{\frac{n(n+1)}{2}}\prod_{i=1}^n\frac{\Gamma
(2p+2i)}{\Gamma (2p+n+i+1)}.
\end{equation}
With  (\ref{JJJ1}), we
have
\begin{equation}\label{lambda1}
\Lambda_1=2^{-n}\pi^{-\frac{n(n+1)}{2}}\prod_{i=1}^n\frac{\Gamma
(k-i)}{\Gamma (k-2i)},
\end{equation}
which is formula (4.3) in \cite{berezin1}.
The holomorphic multiplier representation of the group $\text{Sp}(n, \R )$
on the space of functions defined on the manifold $M=\mc{D}_n$ has the
expression 
$$\pi(g)f(W)=\mu (g,W)f(g^{-1}\cdot W),$$ with the multiplier $\mu$  related
to the cocycle $J$, i.e. $$J(g_1g_2, W)= J(g_1, g_2\cdot W), J(g_2, W), $$ 
 by the relation$$J(g, W)= \mu (g^{-1}, W)^{-1}.$$
For $g\in \text{Sp}( n, \R )$ given by  (\ref{dg}),  the multiplier is 
$$J(g,W)= \det (a^*+ Wb^*)^{\frac{k}{2}}.$$ 
Note that $$K(g\cdot X, g\cdot Y)= J(g,Y)K(X,Y)J(g,X)^*. $$
With  (\ref{dg1}),
(\ref{r4}), we have 
$$g^{-1}\cdot W =
(a^*W-b^t)(-b^*W+a^t)^{-1}=(-W\bar{b}+a)^{-1}(-b+Wa),$$
and finally we have the continuous unitary
 holomorphic  representation of $\text{Sp}(n,\R ) $ on
the space of holomorphic functions of $\mc{D}_n$ attached to the
kernel (\ref{kzz}) (see \cite{god})
\begin{equation}\label{sprep}
\pi \left(\begin{array}{cc}a & b\\ \bar{b}& \bar{a}\end{array}\right)f(W)
= \det (a- W\bar{b})^{-\frac{k}{2}}
f ((-W\bar{b}+a)^{-1}(-b+Wa)), W\in\mc{D}_n.
\end{equation}
Starting from the  development given by Hua 
\cite{hua} (see also \cite{iuri})
   of the
determinant $\det (1-XY^*)^{-\epsilon}$ for the complex Grassmann
manifold in Schur functions \cite{mac},    Berezin \cite{berezin2},
\cite{berezin1} 
found out that
 {\it the admissible set for $k$
for  the
space of functions \fl ~ endowed with the scalar product
\rm{(\ref{ofiW})}}, i.e. the
set of values of $k$ on which the integral, or  it analytic
continuation, converges, for a sufficiently large set of functions
(\ref{ofiW}), 
is the ({\it Wallach} \cite{wallach}) set $\Sigma$
\begin{equation}\label{wal}
\Sigma =  \{ 0, 1, \cdots , n-1\} \cup ( (n-1),\infty ).
\end{equation}
%Note that this set is different of that given in \cite{berezin1} by a 
%factor 2 !
{\it The integral} (\ref{ofiW}) {\it deals with a non-negative scalar
product if $k\ge n-1$, in which the domain of convergence $k \ge 2n$ is
included, and the separate points $k=0,1, ...,n-1$}. The corresponding
coherent states are a ``generalized overcomplete family of states''
(cf. \cite{berezin1}).

In this paper we are not concerned with the analytic continuation of
the discrete series representations. Here are some 
references for  the analytic continuation of the discrete
series or related topics:
\cite{sch,tak,gin,rossi,wallach,up,far1,far,zh,dav}.

\section{The Jacobi group $G^J_n$}\label{JTR}
\subsection{The Jacobi algebra}

 The Jacobi algebra is the  the semi-direct sum
\begin{equation}\label{baza}
\got{g}^J_n:= \got{h}_n\rtimes \got{sp}(n,\R ),
\end{equation}
where $\got{h}_n$ is an  ideal in $\got{g}$,
i.e. $[\got{h}_n,\got{g}]=\got{h}_n$, 
determined by the commutation relations:
\begin{subequations}\label{baza3}
\begin{eqnarray}
\label{baza31}[a^+_k,K^+_{ij}] & = & [a_k,K^-_{ij}]=0,  \\
\label{baza32}[a_i,K^+_{kj}] & = &
\frac{1}{2}\delta_{ik}a^+_j+\frac{1}{2}\delta_{ij}a^+_k 
,\\
 \label{baza33}[K^-_{kj},a^+_i] & = &
\frac{1}{2}\delta_{ik}a_j+\frac{1}{2}\delta_{ij}a_k ,\\
\label{baza34} [K^0_{ij},a^+_k] & = & \frac{1}{2}\delta_{jk}a^+_i,\\
\label{baza35}[a_k,K^0_{ij}] & = & \frac{1}{2}\delta_{ik}a_{j} .
\end{eqnarray}
\end{subequations}

\subsection{Coherent states for the Jacobi group}

 Perelomov's coherent state vectors   associated to the group $G^J_n$ with 
Lie algebra the Jacobi algebra (\ref{baza}), based on the complex
 $N$-dimensional,
$N= \frac{n(n+3)}{2}$, manifold  $M$:
\begin{subequations}\label{nmm}
\begin{eqnarray}\label{mm17}
M& := & H_n/\R\times \text{Sp}(n,\R )/\text{U}(n),\\
\label{mm}
M & = & \mc{D}^J_n:=\C^n\times\mc{D}_n, 
\end{eqnarray}
\end{subequations}
are defined as 
\begin{equation}\label{csu}
e_{z,W}= \exp ({\mb{X}})e_0, 
~\mb{X} := \sum_i z_i a^+_i + \sum_{ij} w_{ij}\mb{K}^+_{ij},~
 z\in \C^n;    W\in\mc{D}_n.
\end{equation} 
The vector $e_0$ verify (\ref{vacuma}) and (\ref{vac}).

\subsection{The differential action}
The differential action of the generators of the Jacobi
group follows from the formulas:
\begin{subequations}
\begin{eqnarray}
 \!\!a^+_ke_{z,W} \!& \!\!\!\!=\!\! & \frac{\pa}{\pa z_k}e_{z,W},\\
\!\!a_ke_{z,W} \!& \!\!\!\!=\!\! & \left( z_k + w_{ki}\frac{\pa}{\pa z_i}\right)e_{z,W},\\
\!\!\mb{K}^0_{kl}e_{z,W}  & \!\!\!\!=\!\! &  \left( \frac{k}{4}\delta_{kl}+
\frac{z_l}{2}\frac{\pa}{\pa z_k} +w_{li}\frac{\pa}{\pa w_{ik}}\right) 
e_{z,W},\\
  \!\!\mb{K}^+_{kl}e_{z,W}  \!& \!\!\!=\!\!\! & \! \!\frac{\pa}{w_{kl}}e_{z,W},\\
 \mb{K}^-_{kl}e_{z,W} \! & \!\!\!\!=\!\! 
 & \left[\frac{k}{2}w_{kl} \!+\! \frac{z_kz_l}{2}
+\frac{1}{2}(z_lw_{ik}\!+\!z_kw_{il})\frac{\pa}{\pa z_i}
\!+\!w_{\alpha l}w_{ki} \frac{\pa}{\pa w_{i \alpha }}\right]\!e_{z,W}, 
\end{eqnarray}
\end{subequations}
and we have
\begin{Remark}\label{difa}
The differential action of the generators of Jacobi group $G^J_n$ is
given by the formulas: 
\begin{subequations}\label{difjac}
\begin{eqnarray}
\mb{a} & = & \frac{\pa}{\pa z},\\
\mb{a}^+ & = & z + W\frac{\pa}{\pa z},\\
\db{K}^- & = & \frac{\pa}{\pa W},\\
\db{K}^0 & = & \frac{k}{4}1 + \frac{1}{2}\frac{\pa}{\pa z}\otimes z 
+\frac{\pa}{\pa W}W,\\
\db{K}^+ & = &  \frac{k}{2}W+\frac{1}{2}z\otimes z 
+\frac{1}{2}(W\frac{\pa}{\pa z}\otimes z+ z\otimes \frac{\pa}{\pa z}W )
+ W\frac{\pa}{\pa W}W.
\end{eqnarray}
\end{subequations}
\end{Remark}
In  (\ref{difjac}) $A\otimes B$ denotes the Kronecker product
 of matrices, here 
$(A\otimes B)_{kl}=a_kb_l$, $A=(a_k), ~B= (b_l), k, l = 1,\dots , n$.

\subsection{Holstein-Primakoff-Bogoliubov type equations}\label{hpb}

We recall the {\it Holstein-Pri\-makoff-Bogoliubov  equations},
%\cite{hol},\cite{bog} (see also \cite{nieto}),
 a consequence
%of the equation (\ref{for}) and 
of the fact that the Heisenberg algebra
is an ideal in the Jacobi algebra %(\ref{baza}), as expressed in
 %(\ref{baza3})-(\ref{baza5}): 
\begin{subequations}\label{hol}
\begin{eqnarray}
\underline{S}^{-1}(Z)\, a_k \, \underline{S}(Z) & = &
 ({\cosh(\sqrt{Z\bar{Z}})\, a})_k + 
(\frac{\sinh (\sqrt{Z\bar{Z}})}{Z\bar{Z}}Z\,a^+)_k,\\
\underline{S}^{-1}(Z)\, a^+ _k\, \underline{S}(Z) & = &
 (\frac{\sinh (\sqrt{\bar{Z}Z})}{\sqrt{\bar{Z}Z}}\bar{Z}a)_k + 
 (\cosh (\sqrt{\bar{Z}Z})\, a^+)_k   ,
\end{eqnarray}
\end{subequations}
and the CCR are still fulfilled in the new creation and annihilation
 operators.
Above $a$ ($a^+$) denotes the column vector formed from
$a_1, ..., a_n$ (respectively,
$a_1^+, ..., a_n^+$).

Let us  introduce the notation: 
\begin{equation}\label{notatie}
\tilde{A}=\left(\begin{array}{c}
A\\ \bar{A}\end{array}\right);~~
\mc{D}={\mc{ D}}(Z) = 
\left( \begin{array}{cc} m & n \\ p & q \end{array}  \right),
\end{equation}
%\begin{equation}\label{notatie1} M  =  \cosh \,(|z|) ;~
%N  =  \frac{z}{|z|}\sinh\,(|z|) ;~
%P  =  \bar{N} ;~
%Q  =   M .
%\end{equation}
\begin{equation}
{\mc{ D}}(Z) = \ee^{X},~\text{where} ~
 X:=\left(\begin{array}{cc} 0 & Z\\ \bar{Z} &
0\end{array}\right)  ,
\end{equation}
where $m,~ n,~ p,~q$ are given by   (\ref{notatie1}).
\begin{Remark}
With the notation {\em (\ref{notatie}), (\ref{notatie1})}, equations
    {\em(\ref{hol})} become:
$$\underline{S}^{-1}(Z)\tilde{a}\underline{S}(Z)={\mc{D}}(Z)\tilde{a}.
$$
\end{Remark}

\begin{Remark}\label{R44} If $D$ and $\underline{S}(Z)$ are defined by
{\rm{(\ref{deplasare})}}, respectively {\rm{(\ref{u1})}}, then 
\begin{equation}\label{schimb0}
D(\alpha )\underline{S}(Z)= \underline{S}(Z) D(\beta ),
\end{equation}
where
\begin{subequations}\label{schimb}
\begin{eqnarray}
\label{schimb33} 
   \beta & = & m\alpha - n\bar{\alpha} =  
   \cosh\,(\sqrt{Z\bar{Z}})\alpha 
-\frac{\sinh(\sqrt{Z\bar{Z}})}{\sqrt{Z\bar{Z}}}Z\bar{\alpha}\\
 & = & (1-W\bar{W})^{-1/2}(\alpha - W \bar{\alpha}),
\end{eqnarray}
\end{subequations}
and
\begin{subequations}\label{schimb34}
\begin{eqnarray}
 \alpha & = & m\beta +  n\bar{\beta} =  
  \cosh\,(\sqrt{Z\bar{Z}})\beta
+\frac{\sinh(\sqrt{Z\bar{Z}})}{\sqrt{Z\bar{Z}}}Z\bar{\beta}\\
 & = & (1-W\bar{W})^{-1/2}(\beta + W \bar{\alpha}).
\end{eqnarray}
\end{subequations}

With the convention {\em (\ref{notatie})}, equation  {\em (\ref{schimb33})}
 can be written down as:
\begin{equation}\label{schimb2}
\tilde{\beta}={\mc{D}}(-Z)\tilde{\alpha};~
 \tilde{\alpha}={\mc{D}}(Z)\tilde{\beta}.
\end{equation}
\end{Remark}

Let us introduce  the notation
\begin{equation}\label{k00}
\underline{S}(Z,A ):= \exp (\sum 2a_{ij}
{\mb{K}}^0_{ij}+z_{ij}{\mb{K}}^+_{ij}-\bar{z}_{ij}{\mb{K}}^-_{ij}). 
\end{equation}
If $g= \left(\begin{array}{ll}\alpha & \beta \\ \bar{\beta} &\bar{\alpha}
\end{array}\right)\in \text{Sp}(n,\R )$,
  then  %%% equations (\ref{hol1}) can be written down as 
\begin{subequations}\label{hol2}
\begin{eqnarray}
S^{-1}(g)\,a\,S(g) & = & \alpha\, a +\beta \, a^+, \\
S^{-1}(g)\,a^+\,S(g) & = & \bar{\beta}\, a +\bar{\alpha} \, a^+ ,
\end{eqnarray}
\end{subequations}
and we have  the following
(generalized Holstein-Primakoff-Bogoliubov) equations:
\begin{Remark}\label{hoho} If $S$ denotes the representation of 
${\rm{{Sp}}}(n,\R )$, with
the convention {\em{(\ref{notatie})}},  we have 
\begin{equation}\label{hol3}
S^{-1}(g)\,\tilde{a}\, S (g)= g\cdot \tilde{a}.
\end{equation}
\end{Remark}

\begin{Remark}\label{remark6}
In the matrix realization of Table 1 % {\em (\ref{nr2})},
%%% equation  {\em(\ref{maimult})} can be written down as
\begin{equation}\label{maimult3}
S (g) D(\alpha ) S^{-1}(g) = D (\alpha_g),
\end{equation}
where %%% {\em (\ref{alpha1})}
one  has the expression of the natural action of
$\rm{Sp}(n,\R )\times \C^n\rightarrow \C^n$: $g\cdot\tilde{\alpha}:=\alpha_g$, 
\begin{equation}\label{alpha2}
\alpha_g = a\,\alpha + b\,\bar{\alpha} .%, 
\end{equation}
% and $a$, $b$ have the expression {\em (\ref{zzz})}.
\end{Remark}

\subsection{The  action of the Jacobi group}\label{actg}

Now we find  a relation  between the normalized
(``squeezed'')  vector
\begin{equation}\label{psi2}
 \Psi_{\alpha, W}:=D(\alpha ) S(W) e_0, ~\alpha\in \C^n, W\in\mc{D}_n,
\end{equation}
and the un-normalized Perelomov's  CS-vector  (\ref{csu}),  
 which is important in the proof of Proposition
\ref{mm1}. %, our   main result of this section.
\begin{lemma}\label{lema6} The vectors 
{\em (\ref{psi2}), (\ref{csu})}, i.e.
$$\Psi_{\alpha, W}:=D(\alpha ) S(W) e_0;~ e_{z,W'}:=\exp
(za^++W'{\mb{K}}^+)e_0 , $$
are related by the relation
\begin{equation}\label{csv}
\Psi_{\alpha, W} = \det (1-W\bar{W})^{k/4}
\exp (-\frac{\bar{\alpha}}{2}z)e_{z,W}  ,
\end{equation}
where 
\begin{equation}\label{cez}
z=\alpha-W\bar{\alpha}.
\end{equation}
\end{lemma}

\begin{proof}
We obtain  successively
\begin{eqnarray*}
\Psi_{\alpha, W}& = &\det (1-W\bar{W})^{k/4}
 D(\alpha ) \exp (W{\mb{K}}^+)e_0\\
 &= &\det (1-W\bar{W})^{k/4}\exp(-\frac{1}{2}|\alpha
|^ 2)
\exp (\alpha a^+)\exp(-\bar{\alpha}a)\exp (W{\mb{K}}^+)e_0\\
 &  = &\det (1-W\bar{W})^{k/4}\exp(-\frac{1}{2}|\alpha
|^2)\exp (\alpha a^+)\exp(-\bar{\alpha}a)\times \\
 & ~ & \exp
(W{\mb{K}}^+)\exp(\bar{\alpha }a)\exp(-\bar{\alpha }a) e_0\\
 & = &\det (1-W\bar{W})^{k/4} \exp(-\frac{1}{2}|\alpha
|^2)\exp (\alpha a^+) E e_0 ,
\end{eqnarray*}
where
\begin{equation}\label{eeee}
E:=\exp(-\bar{\alpha}a)\exp
(W{\mb{K}}^+)\exp(\bar{\alpha }a) .
\end{equation}
As a consequence of  (\ref{bazap}), 
$$\exp (Z)\exp (X) \exp (-Z)=
\exp (X+[Z,X]+\frac{1}{2}[Z,[Z,X]]+\cdots ),$$
where, if we take

$Z=-\bar{\alpha}a;~ X =W{\mb{K}}^+$, then 
$$ [Z,X]=-\bar{\alpha}_k w_{kj}a^+_j= -\bar{\alpha}^tWa^+;$$ 
$$ [Z,[Z,X]]=\bar{\alpha}_pw_{kj}a^+_j=\bar{\alpha}^tW\bar{\alpha} ,$$
where $\alpha^t= (\alpha_1,\dots ,\alpha_n)$. 
We find for $E$ defined by (\ref{eeee}) the value
$$E= \exp (W{\mb{K}}^+-\bar{\alpha}^tWa^++
\frac{\bar{\alpha}^tW\bar{\alpha}}{2}) ,$$
and finally

$$\Psi_{\alpha, W} = 
\exp (-\frac{1}{2}\bar{\alpha}^tz)\det (1-W\bar{W})^{k/4}
e_{\alpha -W\bar{\alpha},W}, $$
i.e.  (\ref{csv}).%%%\hfill $\gata$
\end{proof}
\begin{Remark}\label{comm1}
Starting from {\em  (\ref{csv})}, 
we obtain the expression  {\em (\ref{hot})}
 of the reproducing kernel $K =(e_{z,W},e_{z,W}) $ and {\em{(\ref{KHK})}} of 
$K(z,W;\bar{x},\bar{V}):=
(e_{x,V},e_{z,W})=(e_{\bar{z},\bar{W}},e_{\bar{x},\bar{V}})$. 
\end{Remark}
\begin{proof}
Indeed,   the normalization
$(\Psi_{\alpha ,W},\Psi_{\alpha , W})=1$ imply  that
\begin{equation}\label{are}
(e_{z,W},e_{z,W})=\det(1-W\bar{W})^{-k/2}\exp F,
 F:= \frac{1}{2}(\bar{\alpha}^tz+c.c.).
\end{equation}
With  the notation (\ref{cez}), we have 
$$ \alpha = (1-W\bar{W})^{-1}(z+W\bar{z}), $$
and then $F$ in   (\ref{are}) can be rewritten down as
\begin{equation}\label{are2}
2F = 2\bar{z}^t(1-W\bar{W})^{-1}z+z^t\bar{W}(1-W\bar{W})^{-1}z
+\bar{z}^t(1-W\bar{W})^{-1}W \bar{z},
\end{equation}
\begin{equation}\label{hot}
\!(e_{z,W},e_{z,W}\!)\!=\!\det (M)^{k/2}
\!\exp \frac{1}{2}[2<\!z, M z\!> \!+\! <\!W\bar{z},M z\!>
\!+\!<\!z,M W\bar{z}\!>], 
\end{equation} 
$$ M= (1-W\bar{W})^{-1}.$$
 Above $<x,y>=\bar{x}^ty=\sum\bar{x}_iy_i $. Finally, we find out  
\begin{equation}\label{KHK}
(e_{x,V},e_{y,W}) = \det (U)^{\!k/2}
\exp \frac{1}{2}[2<\!x, U y\!> \!+\! <\!V\bar{y},U y\!>
\!+\!<\!x,U W\bar{x}\!>], 
\end{equation} 
$$ U= (1-W\bar{V})^{-1}.$$
%%%\hfill $\gata$
\end{proof}

From the following proposition we can  see the holomorphic 
action of the  Jacobi group  
\begin{equation}\label{jac}
 G^J_n:=H_n\rtimes \rm{Sp}(n,\R ) 
\end{equation}
on the manifold $M$  (\ref{mm}):

\begin{Proposition}\label{mm1}
Let us consider the action $S(g)D(\alpha )e_{z,W}$, where $g\in
\rm{Sp}(n,\R )$ has the form {\em (\ref{dg})},  $D(\alpha )$ is given by {\em
 (\ref{deplasare})}, and the
coherent state vector is defined in {\em  (\ref{csu})}. Then we have the
formula {\em (\ref{xx})} and the relations {\em (\ref{xxx})-(\ref{x3})}:
\begin{equation}\label{xx}
S(g)D(\alpha )e_{z,W}=\lambda e_{z_1,W_1}, ~ \lambda = \lambda
(g,\alpha; z,W) , \end{equation}
\begin{equation}\label{xxx}
z_1= (Wb^*+ a^*)^{-1}(z+ \alpha -W\bar{\alpha}),
\end{equation}
\begin{equation}\label{x44} W_1=g\cdot
W= (aW+b)(\bar{b}W+\bar{a})^{-1}=(Wb^*+a^*)^{-1}(b^t+Wa^t),\end{equation}
\begin{equation}\label{x1}
\lambda= \det (Wb^*+a^*)^{-k/2}
\exp (\frac{\bar{x}}{2}z
-\frac{\bar{y}}{2}{z_1})
\exp \i\theta_h(\alpha , x),\end{equation}
\begin{equation}\label{x2}x =
(1-W\bar{W})^{-1}(z+W\bar{z}),\end{equation} 
\begin{equation}\label{x3}
y =a (\alpha + x ) + b (\bar{\alpha}+\bar{x}). 
\end{equation}
\end{Proposition}

\begin{proof} 
With  Lemma \ref{lema6}, we have
$e_{z,W}=\lambda_1\Psi_{\alpha_0, W}$, where $\alpha_0$ is given by
 (\ref{x2}) and $\lambda_1=\det (1-W\bar{W})^{-k/4}
\exp (\frac{\bar{\alpha}_0}{2}z)$.
 Then $I:= S(g)D(\alpha)e_{z,W}$ becomes successively
\begin{eqnarray*}
 I & = & \lambda_1 S(g) D(\alpha )\Psi_{\alpha_0,W}\\
 ~ & = & \lambda_1 S(g) D(\alpha ) D(\alpha_0) S(W) e_0\\
~ & = & \lambda_2 S(g) D(\alpha_1 ) S(W) e_0,
\end{eqnarray*}
where $\alpha_1 =\alpha +\alpha_0$ and 
$\lambda_2=\lambda_1\ee^{\i\theta_h(\alpha_1,\alpha_0)}$. With
equations  (\ref{maimult3}),
(\ref{alpha2}), we have $I=\lambda_2 D(\alpha_2)S(g)S(W)e_0$, where
$\alpha_2=a\alpha_1+b\bar{\alpha}_1$. But  (\ref{r1}) implies
$I=\lambda_3 D(\alpha_2)S(g)e_{0,W}$, with
$\lambda_3=\lambda_2\det (1-W\bar{W})^{k/4}$. Now we use  (\ref{r3}) and we
find $I=\lambda_4 D(\alpha_2)e_{0,W_1}$, where, in accord with
 (\ref{r4}), $W_1$ is given by
(\ref{xxx}), and $\lambda_4=\det (Wb^*+ a^*)^{-k/2}\lambda_3$. We
rewrite the last equation as $I= \lambda_5 D(\alpha_2)S(W_1)e_0$,
where $\lambda_5= (1-W_1\bar{W}_1)^{-k/4}\lambda_4$. Then we apply again
Lemma \ref{lema6} and we find $I=\lambda_6e_{z_1,W_1}$, where
$\lambda_6=\lambda_5(1-W_1\bar{W}_1)^k\exp (-\frac{\bar{\alpha}_2}{2}{z_1})$,
and  $z_1=\alpha_2-W_1\bar{\alpha}_2$. Proposition \ref{mm1} is
proved. %%\hfill  $\gata$
\end{proof}

\begin{corollary}The action of the Jacobi group {\em (\ref{jac})} on the
manifold {\em (\ref{mm})}, where $\mc{D}_n=\rm{Sp}(n,\R )/\rm{U}(n)$, is given by
equations {\em (\ref{xx}), (\ref{xxx})}. The composition law in $G$ is
\begin{equation}\label{compositie}
(g_1,\alpha_1,t_1)\circ (g_2,\alpha_2, t_2)= (g_1\circ g_2,
g_2^{-1}\cdot \alpha_1+\alpha_2, t_1+ t_2 +\Im
(g^{-1}_2\cdot\alpha_1\bar{\alpha}_2)),
\end{equation}
where $g\cdot\alpha :=\alpha_g$ is given by {\em  (\ref{alpha2})}, and if
$g$ has the form given by {\em  (\ref{dg})},  then
$g^{-1}\cdot\alpha ={a}^*\alpha -b^t\bar{\alpha}$. 
\end{corollary}
\begin{Remark}\label{re9}
Combining the expressions {\em (\ref{xxx})-(\ref{x3})} and taking
into account the relations {\em (\ref{simplectic})}, the factor
$\lambda$ in {\em{(\ref{xx})}} can be written down as 
\begin{equation}\label{x4}
\lambda = \det (Wb^*+a^*)^{-\frac{k}{2}}\exp(-\lambda_1),
\end{equation}
where 
\begin{equation}\label{x5}
2\lambda_1= z^t(\bar{a}+\bar{b}W)^{-1}\bar{b}z
+(\alpha^t+\bar{\alpha}^t\bar{b}^{-1}\bar{a})(\bar{a}+\bar{b}W)^{-1}\bar{b}
(2z+z_0); z_0=\alpha-W\bar{\alpha},
\end{equation}
or
\begin{equation}\label{x7}
2\lambda_1 = z^t(W+T)^{-1}z+ (\alpha^t+\bar{\alpha}^tT)(W+T)^{-1}(2z+z_0); ~ 
T=\bar{b}^{-1}\bar{a} .
\end{equation}
 In the case $n=1$ the
expression {\em{(\ref{x4})-(\ref{x5})}}
is identical with the expression
 given in {\em Theorem  1.4} in {\em \cite{ez}}
 of the Jacobi forms, under
the the identification of $c, d,\tau, z, \mu, \lambda$ in \cite{ez}
with, respectively,  $\bar{b}, \bar{a}, w, z,\alpha ,-\bar{\alpha}$ in our
notation. Note also that the composition law {\em{(\ref{compositie})}} of the
Jacobi group $G^J$ and the action of the Jacobi group on the base
manifold {\em{(\ref{mm})}} is similar with that 
 in the paper \cite{bb}.  {\em See also %  Remark \ref{rem19} in \S
% \ref{lastc} and
 the Corollary 3.4.4 in \cite{bs}.}
\end{Remark}

\subsection{The K\"ahler  two-form $\omega$ and the volume form}\label{two}
Now we follow the general prescription \cite{last,sinaia}.
We calculate the K\"ahler potential
as the logarithm of the reproducing kernel ({\ref{hot}}),
$f:=\log K$, 
i.e.
\begin{eqnarray}\label{keler}
f  & = & -\frac{k}{2}\log \det (1-W\bar{W})+
\bar{z}_i(1-W\bar{W})^{-1}_{ij}{z}_j+ \\
  & ~ & \nonumber
 \frac{1}{2}[z_i[\bar{W}(1-W\bar{W})^{-1}]_{ij}z_j
+\bar{z}_i[(1-W\bar{W})^{-1}W]_{ij}\bar{z}_j ].
\end{eqnarray}

{\it The K\"ahler two-form $\omega$ on the manifold \text{(\ref{vacuma})} is given by
the formula}:
\begin{equation}\label{aaa}
-\!\i \!\omega \!=\! f_{z_i\bar{z}_j}dz_i\wedge d\bar{z}_j\!+\!
f_{z_i\bar{w}_{\alpha \beta}}dz_i\wedge d\bar{w}_{\alpha \beta}
\!-\!f_{\bar{z}_iw_{\alpha \beta}}d\bar{z}_i\wedge d w_{\alpha\beta}
 \!+\!f_{w_{\alpha\beta}\bar{w}_{\gamma \delta}}dw_{\alpha\beta}\wedge 
d\bar{w}_{\gamma \delta},
\end{equation}
{\it where}
\begin{subequations}\label{aas}
\begin{eqnarray}
f_{z_i\bar{z}_j} & = & (1-W\bar{W})^{-1}_{ji}\label{aa1}, \\
2f_{\bar{z}_i{w}_{\alpha \beta} }& = &
  2(1-W\bar{W})^{-1}_{i\alpha}[\bar{W}(1-\bar{W}{W})^{-1}z]_{\beta}
\label{aa2} \\
 & & + (1-W\bar{W})^{-1}_{i\alpha}[(1-W\bar{W})^{-1}\bar{z}]_{\beta}
\nonumber\\
& &  + [\bar{z}^t(1-W\bar{W})^{-1}]_{\alpha}(1-\bar{W}W)^{-1}_{\beta i},
\nonumber\\
2 f_{w_{\alpha \beta}\bar{w}_{\gamma\delta}} & = &
k (1-\bar{W}W)^{-1}_{\beta\gamma}(1-W\bar{W})^{-1}_{\delta\alpha}
\label{aa3}\\
& & + 2[\bar{z}^t(1-W\bar{W})^{-1}W]_{\gamma}
(1-W\bar{W})^{-1}_{\delta\alpha}[\bar{W}(1-W\bar{W})^{-1}z]_{\beta}\nonumber
\\
 & &  +
2[\bar{z}^{t}(1-W\bar{W})^{-1}]_{\alpha}(1-\bar{W}W)^{-1}_{\beta\gamma}
[(1-W\bar{W})^{-1}z]_{\delta}\nonumber\\
& & +  [z^t(1-W\bar{W})^{-1}]_{\gamma}(1-W\bar{W})^{-1}_{\delta\alpha}
[(1-\bar{W}W)^{-1}z]_{\beta}\nonumber\\
& & +[z^t(1-W\bar{W})^{-1}]_{\alpha}(1-W\bar{W})^{-1}_{\delta\beta}
[z^t(1-W\bar{W})^{-1}W]_{\gamma}\nonumber \\
 & & + [\bar{z}^t(1-W\bar{W})^{-1}W]_{\gamma}
(1-W\bar{W})^{-1}_{\delta\alpha}[(1-\bar{W}W)^{-1}\bar{z}]_{\beta}\nonumber\\
 & & + [\bar{z}^t(1-W\bar{W})^{-1}]_{\alpha}(1-\bar{W}W)^{-1}]_{\beta\gamma}
[W(1-\bar{W}W)\bar{z}]_{\delta},\nonumber
\end{eqnarray}
\end{subequations}
i.e.
\begin{eqnarray}\label{aaa1}
-\i \omega  & = & \frac{k}{2}\tr [(1-W\bar{W})^{-1}d W
\wedge (1-\bar{W}W)^{-1} d\bar{W}]\\
 & & + ~ \tr [dz^t\wedge (1-\bar{W}W)^{-1} d\bar{z}] \nonumber\\
 & &  - ~ \tr [d\bar{z}^t(1-W\bar{W})^{-1}\wedge dW \bar{x}] + c.c. \nonumber\\
& &  +  ~ \tr [\bar{x}^tdW(1-\bar{W}W)^{-1}\wedge d\bar{W} {x}]\nonumber
\end{eqnarray}
where $x$ is defined in  (\ref{x2}).

We can write down the two-form $\omega$ (\ref{aaa1}) as
\begin{equation}\label {aab}
-\!\i\omega \!=\! \frac{k}{2}\tr (B\! \wedge\! \bar{B})
 \!+\!\tr (A^t\bar{M}\!\wedge \!\bar{A}), ~ A\!=\!dz\!+\!dW\bar{x},
~B\!=\! MdW, ~ M\!=\!(1\!-\!W\bar{W})^{-1}. 
\end{equation}

Now we determine the Liouville  form. We apply the following  technique 
(see Hua's book, Ch. IV). {\it  Let $z'=f (g,z)$ the action of the group
$G$ on the circular domain $M$. Let us determine the element $g\in G$ 
such that $z'(z_1)=0$ Then the density of the volume form is
$Q= |J|^2$, where $J$ is the Jacobian $J=\pa z'/\pa z$}.
 We apply this method
to our manifold (\ref{nmm}) and the Jacobi group $G^J_n$ (\ref{jac}). 

The transformation with the desired properties is:
\begin{subequations}\label{transform}
\begin{eqnarray}\label{ttr}
z' & = & U(1-W_1\bar{W}_1)^{1/2}(1-W\bar{W}_1)^{-1}\times\label{trans1}\\
 & &  [z-(1-W\bar{W}_1)(1-W_1\bar{W_1})^{-1}z_1 + 
(W-W_1)(1-\bar{W}_1W_1)^{-1}\bar{z}_1],\nonumber\\
 W' & = & U (1\!-\!W_1\bar{W}_1)^{-1/2}(W\!-\!W_1)
(1\!-\!\bar{W}_1W)^{-1}(1\!-\!\bar{W}_1W_1)^{1/2}U^t\label{trans2},
\end{eqnarray}
\end{subequations}
where $U$ is a unitary matrix. 

We find  that
\begin{equation}\label{tt1}
\pa z'/\pa z=  U(1-W_1\bar{W}_1)^{1/2}(1-W\bar{W}_1)^{-1},
\end{equation}
\begin{equation}\label{tt2}
dW'=A dW A^t, ~ A = U(1-W_1\bar{W}_1)^{1/2}(1-W\bar{W}_1)^{-1}. 
\end{equation}

In order to calculate the Jacobian of the transformation (\ref{tt2}), 
we use the following property extracted from p. 398 in Berezin's paper
\cite{berezin2}:
{\it Let $A$ be a matrix and $L_A$ the transformation of  a matrix of the same order $n$, 
$L_A\xi = A\xi A^t$. If the matrices $A$ and $\xi$ are symmetric, then 
$\det L_A= (\det A)^{n+1}$}. 
The overall determinant of the transformation (\ref{transform}) is
\begin{equation}\label{ddt}
J = \left|\begin{array}{cc}\pa z'/\pa z &\pa z'/\pa W   \\ \pa W'/\pa z
& \pa W'/\pa W \end{array}\right|= \pa z'/\pa z \pa W'/ \pa W,
\end{equation}
because $\pa W'/\pa z =0$. 
Finally, taking $W_1= W$, we find out  
\begin{equation}\label{QQQ}
 Q = \det (1-W\bar{W})^{-(n+2)}.
\end{equation}

\subsection{The scalar product}

If $f_{\psi}(z)= (e_{\bar{z}},\psi )$, then
\begin{equation}\label{ofi}
(\phi ,\psi )= \Lambda\! \int_{z\in\C^n;
1-W\bar{W}>0}\!\bar{f}_{\phi}(z,W)f_{\psi}(z,W)Q K^{-1} dz dW.
\end{equation}
$Q$ is the density of the volume form given by  (\ref{QQQ}), $K$ is
 the reproducing kernel (\ref{hot}), and
\begin{equation}\label{masurae}
dz = \prod_{i=1}^n \Re z_i\Im z_i; ~~ dW = \prod_{1\le i\le j \le n}
\Re w_{ij}\Im {w}_{ij}.
\end{equation}

We have $K^{-1}=\det (1-W\bar{W})^{k/2}\exp (-F)$ with $F$ given by 
 (\ref{are2}). 

In order to find the value of the constant $\Lambda$ in
 (\ref{ofi}), we take the functions $\phi,\psi=1$, we change the variable
$z =  (1-W\bar{W})^{1/2}x$ and we get
$$1 \!= \!\Lambda\! \int_{1-W\bar{W}>0}\!\!\det (1-W\bar{W})^{\frac{k-3}{2}-n}d W
\!\int_{x\in\C^n}\!\!\exp({-|x|^2})
\exp(-\frac{x^t\!\cdot\!\bar{W}x +\bar{x}^t\!\cdot\! W\bar{x}}{2}) d x . $$

We apply equations (A1), (A2) in  Bargmann \cite{bar70}:
$$I(B,C)=\int \exp
(\frac{1}{2}(x.Bx+\bar{x}.C\bar{x}))\pi^{-n}e^{-|x|^2}
\prod_{k=1}^n d\Re x_k d\Im x_k= [\det (1-CB)]^{-\frac{1}{2}}, $$
where $B$, $C$ are complex symmetric matrices such that $|B|<1,
|C|<1$. Here  $B=-\bar{W}$, $C= -W$. So, we get
$$1\!=  \!\Lambda\pi^n \!\int_{ 1-W\bar{W}>0}\det (1-W\bar{W})^{p}d W , ~ p=
\frac{ k-3}{2} -n.$$
We apply Theorem 2.3.1 p. 46 in Hua's book \cite{hua}
$$\int_{1-W\bar{W}>0, W= W^t} \!\det (1-W\bar{W})^{\lambda} dW=
J_n(\lambda), $$
 and we find out:
\begin{equation}\label{ofi2}
\Lambda = \pi^{-n}J^{-1}_n(p), 
 \end{equation}
where, for $p>-1$, $J_n(p)$ is given by formula (\ref{JJJ}) or by
formula (\ref{JJJ1}).  
So, we find out for $\Lambda$
\begin{equation}\label{ofi1}
\Lambda = \frac{k-3}{2\pi^{\frac{n(n+3)}{2}}}\prod_{i=1}^{n-1}
\frac{(\frac{k-3}{2}-n+i)\Gamma (k+i-2)}{\Gamma [k+2(i-n-1)]} .
\end{equation}

\begin{Proposition}\label{final}
Let us consider the Jacobi group $G^J_n$ 
 {\em (\ref{jac})}  with the composition
rule
{\em (\ref{compositie})}, acting on the
coherent state manifold {\em (\ref{nmm})} via equation 
 {\em  (\ref{xxx})}, {\em  (\ref{x44})}. The manifold
$\mc{D}^J_n$ has the K\"ahler potential {\em (\ref{keler})}  and the
$G^J_n$-invariant K\"ahler two-form $\omega$ given by 
 {\em (\ref{aab})}.
% The holomorphic polynomials {\em(\ref{x3x1})} associated to the coherent
%state vectors {\em (\ref{csu})} are given by {\em  (\ref{x4x})}, where the
%functions
% $f$ are
%given by  {\em (\ref{f2})}, while the polynomials $P$ are given by
%{\em  (\ref{marea})}.
 The Hilbert space of holomorphic functions $\mc{F}_K$
associated to the holomorphic kernel $K:M\times \bar{M}\rightarrow\C$ given
 by 
 {\em (\ref{KHK})} is endowed with the
scalar product {\em (\ref{ofi})}, where the normalization
 constant $\Lambda$ is given by
{\em (\ref{ofi1})} and the density of 
volume given by {\em{(\ref{JJJ})}}.  % $G$-invariant measure $d\nu$
  % {\em (\ref{ofi3})}.
\end{Proposition}

Recalling Proposition IV.1.9. p. 104 in \cite{neeb}, Proposition
\ref{mm1} can be formulated as follows:
\begin{Proposition}\label{finalf}
Let $h := (g,\alpha )\in G^J_n$,  where $G^J_n$ is the Jacobi group
{\em (\ref{jac})}, and we consider the representation
  $\pi (h) := S(g)D(\alpha )$, $g\in
\rm{Sp}(n, \R )$,
$\alpha\in \C^n$, and let the notation $x:=(z,W)\in
{\mc{D}}^J_n:=\C^n\times\mc{D}_n$. Then the 
continuous unitary representation $(\pi_K,\Hi_K)$  attached to the
positive definite holomorphic kernel $K$ defined by {\em  (\ref{KHK})} is
\begin{equation}\label{rep}
(\pi_K(h).f)(x)=J(h^{-1},x)^{-1}f(h^{-1}.x),
\end{equation}
where the cocycle $J(h^{-1},x)^{-1}:=\lambda (h^{-1},x)$  with
$\lambda$ defined by
{\em  (\ref{xx})-(\ref{x3})} and the function $f$
belongs to the Hilbert space of holomorphic functions $\Hi_K\equiv
\mc{F}_K$ endowed with the
scalar product {\em (\ref{ofi})}, where $\Lambda$ is given by 
{\em  (\ref{ofi1})}.
\end{Proposition}
\begin{Comment}
The value of $\Lambda$ {\em{(\ref{ofi})}} given by {\em{(\ref{ofi1})}}
  corresponds to the one given in equation {\em{(7.16)}} in \cite{holl}, taking above
$ n=1$, $k\rightarrow 4k$.  Note that $p$ defining
 the normalization
constant $\Lambda$  in  {\em{(\ref{ofi2})}} for the Jacobi group
 is related with  $q$  in  {\em{(\ref{ofiW})}} 
defining the normalization
 constant for the group $\rm{Sp}(n,\R )$ by the relation $p =q-\frac{1}{2}$. 
\end{Comment}

\section{The Jacobi group $G^J_1$}\label{unu1}

\subsection{K\"ahler-Berndt's K\"ahler two-form $\omega$}\label{abcd}
We recall \cite{holl} that in the case $n=1$ formula (\ref{aab}) reads:
\begin{equation}\label{aab1}
-\i\omega =\frac{2k}{(1-w\bar{w})^2}dw \wedge d\bar{w} +
\frac{A\wedge \bar{A}}{1-w\bar{w}},
~A=dz+\bar{\alpha}_0dw, ~\alpha_0=\frac{z+\bar{z}w}{1-w\bar{w}}.
\end{equation}

 Rolf Berndt -alone or in collaboration - has studied the 
 real Jacobi group $G^J(\R )_1$ in
several references, from which I mention 
\cite{bern,bb,bs}. The Jacobi group appears (see explanation in
\cite{cal}) in the
context of the so called {\it Poincar\'e group} or {\it The New
Poincar\'e group} investigated by Erich K\"ahler as the 10-dimensional
group $G^K$ which invariates a hyperbolic metric
\cite{cal1,cal2,cal3}. K\"ahler and Berndt have investigated the
Jacobi group  $G^J_0(\R ):= \text{SL}_2(\R )\ltimes \R^2$ acting on the
manifold $\mc{X}^J_1:=\mc{H}_1\times\C$, where  $\mc{H}_1$ is the upper half plane
$\mc{H}_1:=\{v\in\C|\Im (v)>0\}$. We recall that
%\begin{Remark}\label{actnea}
{\it the action of $G^J_0(\R )$ on $\mc{X}^J_1$ is
given by $(g,(v,z))\rightarrow (v_1,z_1)$, $g=(M , l )$, where}
\begin{equation}\label{lac}
v_1=\frac{av+b}{cv+d},  z_1=\frac{z+l_1v+l_2}{cv+d}
; ~M  = \left(\begin{array}{cc} a & b\\c & d
\end{array}\right)\in \rm{SL}_2(\R) , (l_1,l_2)\in\R^2 .
\end{equation} 
%\end{Remark}
It can be proved  \cite{holl} that
\begin{Remark}
When expressed in the  coordinates $(v, u)\in\mc{X}^J_1$
 which are related to the
coordinates $(w,z)\in\mc{D}^J_1$ by the map  
\begin{equation}\label{noiec}
w=\frac{v-\i}{v+\i};~ z=\frac{2\i u}{v+\i},
w\in\mc{D}_1,~v\in\mc{H}_1,z\in\C ,
\end{equation}
the  K\"ahler two-form {\em{(\ref{aab1})}} is 
identical with the one  
\begin{equation}\label{ura3}
-\i\omega = \!
-\frac{2k}{(\bar{v}-v)^2}dv\!\wedge\! d\bar{v}
+\!\frac{2}{i(\bar{v}\!-\!v)}B\!\wedge\! \bar{B}, 
~B=
du-\!\frac{u\!-\!\bar{u}}{v\!-\!\bar{v}}dv,
v,u\in\C, \Im (v)\!>0,
\end{equation}
 considered by K\"ahler-Berndt \cite{cal1,cal2,cal3,cal,bern,bs}.
If we use the EZ \cite{bs,ez} coordinates 
adapted to our notation
\begin{equation}\label{eqez}
v=x+\i y; ~ u=pv+q,~ x,p,q, y\in\R ,y>0,
\end{equation}
the $G^J_0(\R )$-invariant K\"ahler metric on $\mc{X}^J_1$
 corresponding to the K\"ahler-Berndt's
K\"ahler two-form $\omega$ {\rm{(\ref{ura3})}} reads
\begin{equation}\label{victorie}
ds^2=\frac{k}{2y^2}(dx^2+dy^2)+\frac{1}{y}[(x^2+y^2)dp^2+dq^2+2xdpdq],
\end{equation}
i.e. the equation at {\em{p. 62}} in \cite{bs} or the equation
given  at {\rm{p. 30}} in \cite{bern}.
\end{Remark}
\subsection{The base of functions}\label{basef}
We reproduce the following result proved in \cite{holl}
\begin{Remark}
The reproducing kernel $K:\mc{D}^J_1 \times \bar{\mc{D}}^J_1 \rightarrow\C$
 has the property:
%%%%\begin{subequations}\label{ker3}
\begin{eqnarray*}
\label{sub1}K(z,w;\bar{z},\bar{w}')
 & := & (e_{\bar{z},\bar{w}},e_{\bar{z}',\bar{w}'})
= \sum_{n,m}f_{|n>,e_{k',k'+m}}(z,w)\bar{f}_{|n>,e_{k',k'+m}}(z',w') \\
\label{sub2} & = & (1-{w}\bar{w}')^{-2k}\exp{\frac{2\bar{z}'{z}+z^2\bar{w}'
+\bar{z}'^2w}{2(1-{w}\bar{w}')}} .
\end{eqnarray*}
%%%\end{subequations}
Here 
\begin{equation}\label{x4x}
f_{{|n>;e_{k',k'+m}}}(z,w)=f_{e_{k',k'+m}}(w)\frac{P_n(z,w)}{\sqrt{n!}}, 
\end{equation}
\begin{equation}\label{f2}
f_{e_{k,k+n}}(z):=(e_{\bar{z}},e_{k,k+n})
=\sqrt{\frac{\Gamma (n+2k)}{n!\Gamma (2k)}}z^n ,
\end{equation}
\begin{equation}\label{x6x}
P_n(z,w) := (\frac{\i}{\sqrt{2}})^n w^{\frac{n}{2}}H_n(\frac{-\i z}{\sqrt{2w}}),
\end{equation} i.e.  the polynomials $P_n(z,w)$ have the expression
\begin{equation}\label{marea}
P_n(z,w)=n!\sum _{k=0}^{[\frac{n}{2}]}
(\frac{w}{2})^k\frac{z^{n-2k}}{k!(n-2k)!} .
\end{equation} 

The first 6
polynomials $P_n(z,w)$ are:
\begin{equation}
\begin{array}{ll}
P_0= 1; & P_1 = z ;\\
P_2 = z^2 +w ; & P_3 = z^3 + 3 zw ; \\
P_4= z^4 + 6 z^2w +3 w^2 ; & P_5 = z^5 +10 z^3 w + 15 z w .
\end{array}
\end{equation}
\end{Remark}

\subsection{Digression on  Jacobi and  Schr\"odinger groups}

$\bullet$ We have called the algebra (\ref{baza}), the {\it Jacobi
algebra} and the
group (\ref{jac}), the {\it Jacobi group}, in agreement  with the name
used   in \cite{bs} or at p. 178 in \cite{neeb}, 
where the algebra $\got{g}^J_1:= 
\got{h}_1\rtimes\got{sl}(2,\R )$  is called ``Jacobi  algebra''.
 The denomination adopted in the present paper is of course
in accord with the one used in \cite{neeb}   because of the
isomorphism of the Lie algebras $\got{su}(1,1)
\sim\got{sl}(2,\R)\sim\got{sp}(1,\R )$$(\sim \got{so}(2,1))$. 
 Also the name ``Jacobi
algebra'' is used in \cite{neeb} p. 248 to call the semi-direct sum of
the $(2n+1)$-dimensional Heisenberg algebra and the symplectic algebra,
 $\got{hsp}:=
\got{h}_n\rtimes \got{sp}(n,\R )$. 
 The group corresponding to this algebra
 is called sometimes  in 
the  Mathematical Physics literature (see e. g. \S 10.1 in \cite{ali},
which is based on \cite{si})
the   {\it metaplectic group}, but in reference \cite{neeb} the term 
``metaplectic group''  is reserved to  the 2-fold covering group of the
symplectic group, cf. p. 402 in  \cite{neeb} (see also \cite{bar70}
and \cite{itzik}). Other names of the metaplectic representation are
the {\it  oscillator representation, the harmonic representation} or
{\it the Segal-Shale-Weil representation}, see references in Chapter 4 of
\cite{fol} and \cite{bs}.

$\bullet$
 Apparently, the denomination of the semidirect product of the Heisenberg
and symplectic group as Jacobi group was introduced by Eichler and Zagier,
cf. Chapter I in \cite{ez}. In  their monograph \cite{ez}, Eichler
and Zagier
 have
introduced the notion of {\it Jacobi form} on $\text{SL}_2(\Z )$ as a
holomorphic function on 
$\mc{X}^J_1$  satisfying three properties. One of this properties,
generalized to other groups, was studied by Pyatetskii-Shapiro, who
referred to it as the {\it Fourier-Jacobi} 
 expansion (see \S 15 in \cite{ps}), and to some
coefficients as {\it Jacobi forms}, a name adopted by Eichler and
Zagier to denote also the group appearing in this context. The
denomination Jacobi group 
was adopted also in the monograph \cite{bs} and this group is
important in K\"ahler's approach, see Chapter 36 in \cite{cal3}.  

$\bullet$ The Jacobi algebra, denoted $ \got{st}( n,\R )$
 by Kirillov
 in \S 18.4 of \cite{kir} or $\got{tsp}(2n+2,\R )$ in \cite{kir2},
 is isomorphic with the subalgebra of Weyl algebra  $A_n$ 
 of polynomials of degree maximum 2 in the variables
$p_1,\dots,p_n, q_1,\dots,q_n$ with the Poisson bracket, while  the
 Heisenberg algebra $\got{h}_n$ is the nilpotent ideal isomorphic with
polynomials 
of degree  $\le 1$ and  the real symplectic algebra $\got{sp}(n,\R
)$ is isomorphic to the subspace of symmetric homogeneous polynomials
of degree $2$.

$\bullet$
  U. Niederer 
 has introduced in \cite{nied} and developed in
\cite{nied2} the  concept of {\it the maximal kinematical invariance group}
 (MKI) of the (free)
Schr\"odinger equation. The  Schr\"odinger group is a 12-parameter group
containing, in addition to the Galilei group $G^G_3$,
the group of dilations
and a 1-parameter group of transformations, called expansions, which is,
in some respects, very similar to the special conformal transformations
of the conformal group.
Let us consider the wave equation 
\begin{equation}\label{nie}
\Delta (t,\mb{x})\psi (t,\mb{x})= 0,
\end{equation} where $\Delta (t,\mb{x})$ is
a wave operator, in particular the Schr\"odinger operator 
$$\Delta
(t,\mb{x})=\i\partial_0+\frac{1}{2m}\Delta_3,~(t,\mb{x})\in\R^4, $$
and let 
\begin{equation}\label{marebranza}(t,\mb{x})\rightarrow g(t,\mb{x})
\end{equation} be any invertible
coordinate transformation $g$. Equation (\ref{nie})
is invariant under the transformation (\ref{marebranza}) if
$$\psi (t,\mb{x})\rightarrow 
(T_g\psi )(t,\mb{x})=f_g[g^{-1}(t,\mb{x})]\psi[g^{-1} (t,\mb{x})],$$
is again a solution of the wave equation when  $\psi$ is a 
wave function of (\ref{nie}).    
Niederer has determined the (one-cocycle) $f_g$ and the
Schr\"odinger transformation \begin{equation}\label{gali}
g(t,\!\mb{x}\!)\!=\!\left(\!d^2\frac{t\!+\!b}{1\!+\!\alpha
(t\!+\!b)},d\frac{R\mb{x}\!+\!\mb{v}t\!+\!\mb{a}}{1\!+\!\alpha
(t\!+\!b)}\!\right) , \alpha,b, d\in\R, R\in \!\text{SO}(3),\!
\mb{v},\!\mb{a}\in\!\R^3. 
\end{equation} 
In fact, $T_g$ is a projective representation of the Schr\"odinger
group. For the Schr\"odinger group  see also \cite{barut}.
Another names for the Schr\"odinger group are the {\it Hagen group} or the
{\it conformal Galilean group} \cite{hagen}. Other references on the
same subject are 
\cite{burdet,roman,barut2,peroud,burdet2,barut3,fein}.

In \cite{dobre} is given the  Levi-Malcev decomposition of the
 Schr\"odinger group in $n+1$-space-time dimensions:
$$\text{Sch}(n)=\left[\underbrace{\R^n\times G^G_n}_{\text{radical}}\right]
\rtimes\left[\underbrace{\text{SL}(2,\R)
\times \text{SO}(n)}_{\text{SS-Levi part}}\right],$$
where the Galilei group is 
$$G^G_n=(\text{SO}(n)\rtimes\R^n)\rtimes\R^{n+1} .$$
In the case $n=1$, (\ref{gali})
with {\bf $t\in\C, \Im (t)>0$}, {\bf $x\in\C$}, ($d=1, R=1$)
corresponds to (\ref{lac}) with 
$$ M  = \left(\begin{array}{cc} 1 & b\\ \alpha & 1+\alpha b,
\end{array}\right)\in \rm{SL}_2(\R) , (l_1,l_2)=(v,a)\in\R^2 .$$
It it interesting that the  paper \cite{cal3} starts with the Chapter
{\it Relativit\"at nach Galilei}. The relationship between the Jacobi
group $G^J$, the New Poincar\'e group $G^K$, the de Sitter group
$\text{SO}_0(1,4)$ and the anti de
Sitter group $\text{SO}_0(2,3)$, the
 (standard) Poincar\'e group  $G^P = \text{SO}_0(1,3)\ltimes\R^4$ and his 
contraction, the
Galilei group $G^P_3$, is discussed at p. 33 in \cite{cal}.

\section*{Acknowledgments}
The author is grateful to 
 Tudor Ratiu for
the hospitality  at the Bernoulli Center, EPFL Lausanne,
 Switzerland, where this investigation was
started in November -- December 2004.   
 The author is thankful to
 the organizers of 
 {\it XXIII Workshop  on Geometric methods in Physics},
  June 26 -- July 2, 2005,   Bia\l owie\.{z}a,  
Poland  and of the  {\it 3$^{~\text{rd}}$ Operator Algebras
and Mathematical Physics Conference}, Bucharest, Romania, 
 August 10 -- 17,  2005,
  for the opportunity to
report results on this subject. 

%\today

\end{document}